\documentclass[a4paper, 12pt]{article}
\usepackage{amssymb}
\usepackage{amsmath, amsthm, amssymb}
\usepackage{amssymb,amsmath,amsthm,url}
\usepackage{enumerate}
\usepackage{verbatim}

\usepackage[all]{xy}
\reversemarginpar
\def\Hom{\mbox{Hom}}

\def\H{\mbox{H}}

\def\fA{\mathfrak{A}}
\def\arcs{\fA}
\def\fU{\mathfrak{U}}
\def\arcsU{\fU}
\def\fd{\mathfrak{d}}
\def\arcsd{\fd}

\def\fX{\mathfrak{X}}
\def\arcsX{\fX}
\def\fY{\mathfrak{Y}}
\def\arcsY{\fY}
\def\fW{\mathfrak{W}}
\def\arcsW{\fW}

\theoremstyle{plain}
 \newtheorem{thm}{Theorem}[section]
 \numberwithin{equation}{section} 
 \numberwithin{figure}{section} 

 \theoremstyle{plain}
 \newtheorem{lem}[thm]{Lemma} 
 \newtheorem{corol}[thm]{Corollary} 
 
 \newtheorem*{mythm}{Theorem}

 \theoremstyle{definition}
 \newtheorem{defn}[thm]{Definition}
 \newtheorem{property}[thm]{Property}

 \theoremstyle{definition}

 \theoremstyle{remark}

 \newtheorem{remark}[thm]{Remark}

 \newtheorem*{acknowledgement}{Acknowledgement}

\title{A characterization of torsion theories in the cluster category of Dynkin type $A_{\infty}$}
\author{ Puiman Ng }
\date{ }

\begin{document}


\setlength{\parindent}{0pt}
\setlength{\parskip}{7pt}

\setcounter{section}{-1}

\maketitle

\pagenumbering{roman}


\pagenumbering{arabic}

\begin{abstract}
Let $\cal D$ be the cluster category of Dynkin type $A_{\infty}$. 
Following the line of \cite{PJ11}, this paper provides a bijection between torsion theories in $\cal D$ and certain configurations of arcs connecting non-neighbouring integers.



\end{abstract}


\section{Introduction}
The cluster category $\cal D$ of Dynkin type $A_{\infty}$ was introduced in \cite{PJ11}. One of its several definitions, which is completely analogous to the definition of the cluster category of type $A_n$, motivates us to say that $\cal D$ is a cluster category of type $A_{\infty}$. Namely, it is the orbit category $D^{f}(\mbox{mod}\,\Gamma)/S\Sigma^{-2}$. Here $\Gamma$ is a quiver of type $A_{\infty}$ with zigzag orientation and $S$ and $\Sigma$ are the Serre and suspension functors of the finite derived category $D^{f}(\mbox{mod}\,\Gamma)$. 

There are also several other ways to realize the category $\cal D$. In brief, it is the algebraic triangulated category generated by a $2$-spherical object. It is also the compact derived category $D^{c}(A)$ of the differential graded cochain algebra $A = C^*(S^2;k)$ where $S^2$ is the $2$-sphere and $k$ is a field. Finally, $\cal D$ is the finite derived category $D^{f}(k[T])$ where $k[T]$ is viewed as a DG algebra with $T$ placed in homological degree $1$ and zero differential. It is ubiquitous and the reader can refer to \cite[Section 0]{PJ11} for more details.

In \cite{PJ11}, the cluster tilting subcategories of $\cal D$ were shown to be in bijection with certain maximal sets of non-crossing arcs connecting non-neighbouring integers. One can think of these maximal sets as ``triangulations of the $\infty$-gon''. 

A torsion theory in $\cal D$ is a pair $(\cal X, \cal Y)$ of subcategories such that there are no non-zero morphisms from any object in $\cal X$ to any object in $\cal Y$, and that for each $d$ in $\cal D$, there is a distinguished triangle $x \rightarrow d \rightarrow y \rightarrow$ in $\cal D$, with $x$ in $\cal X$ and $y$ in $\cal Y$. Note that if $\cal T$ is a cluster tilting subcategory, then ($\cal T$, $\Sigma \cal T$) is a torsion theory, but t-structures and co-t-structures are other examples of torsion theories. 

In this paper, we will generalize the results of \cite{PJ11} by giving a bijection between torsion theories in $\cal D$ and certain configurations of arcs connecting non-neighbouring integers.

Let us recall some material from \cite{PJ11}. The category $\cal D$ has finite dimensional $\Hom$ spaces over a field $k$ and
split idempotents, so it is Krull-Schmidt (\cite[Remark 1.2]{PJ11}). Throughout the paper, a subcategory of $\cal D$ is assumed to be a full subcategory closed under direct sums and direct summands. For subcategories $\cal X$ and $\cal Y$ of $\cal D$, the set of morphisms from any $x$ in $\cal X$ to any $y$ in $\cal Y$ is denoted by $(\cal X,\cal Y)$. In particular, for objects $x$ and $y$ of $\cal D$, the set of morphisms from $x$ to $y$ is denoted by $(x,y)$. 



By \cite[Remark 1.4]{PJ11}, the Auslander-Reiten (AR) quiver of $\cal D$ is $\mathbb{Z} A_{\infty}$, and we will use the following standard coordinate system on the AR quiver.
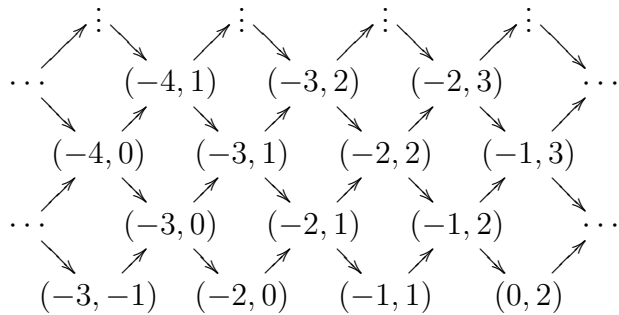
\[
  \xymatrix @-4.0pc @! {
    &  & & \vdots \ar[dr] & & \vdots \ar[dr] & & \vdots \ar[dr] & & \vdots \ar[dr] & &  & \\
     & & \cdots  \ar[ur] \ar[dr] & & (-4,1) \ar[ur] \ar[dr] & & (-3,2) \ar[ur] \ar[dr] & & (-2,3) \ar[ur] \ar[dr] & &  \cdots & &  \\
    &   & & (-4,0) \ar[ur] \ar[dr] & & (-3,1) \ar[ur] \ar[dr] & & (-2,2) \ar[ur] \ar[dr] & & (-1,3) \ar[ur] \ar[dr] & &  & \\
    & &  \cdots \ar[ur]\ar[dr] & & (-3,0) \ar[ur] \ar[dr] & & (-2,1) \ar[ur] \ar[dr] & & (-1,2) \ar[ur] \ar[dr] & & \cdots  & & \\\
    & & & (-3,-1) \ar[ur] & & (-2,0) \ar[ur] & & (-1,1) \ar[ur] & & (0,2) \ar[ur] & &  & \\
               }
\]
Let $\Sigma$ be the suspension functor of $\cal D$. Since $\cal D$ is $2$-Calabi-Yau, its Serre functor is $S = \Sigma^2$ and the AR translation is $\tau = S\Sigma^{-1} = \Sigma$. Note that in terms of coordinates, the action of $\Sigma = \tau$ is $\Sigma(m,n) = (m-1,n-1)$, see \cite[Remark 1.4]{PJ11}.

Now let us have a few definitions to interpret these coordinate pairs in an alternative way, as arcs connecting non-neighbouring integers.

\begin{defn} (\cite[Definition 3.1]{PJ11})
An arc is a pair $(m,n)$ of integers with $n-m \geq 2$.
The arc $(m,n)$ is said to end in each of the integers $m$ and $n$.
Two arcs $(m_1,n_1)$ and $(m_2,n_2)$ are said to cross if we have either $m_1 < m_2 < n_1 < n_2$ or $m_2 < m_1 < n_2 < n_1$. Note that the action of $\Sigma$ makes sense on arcs as well. 
\end{defn}

In all the diagrams to follow arcs are drawn on number lines which are numbered thus:
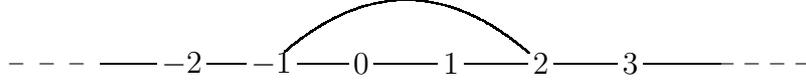
\[
  \xymatrix @-2.5pc @! {
    \rule{0ex}{6.5ex}\ar@{--}[r]
   &\ar@{-}[r]& *{\, -2 \, } \ar@{-}[r]
& *{\, -1 \, } *{} \ar@{-}[r] \ar@/^2pc/@{-}[rrr]
& *{\, 0 \, } \ar@{-}[r] & *{\, 1 \, } \ar@{-}[r] 
& *{\, 2 \, } \ar@{-}[r] & *{\, 3 \, } \ar@{-}[r] &*{}\ar@{--}[r]&*{}
                    }
\]

For example, in the above diagram, the arc $(-1,2)$ is drawn as a curve between the integers $-1$ and $2$. Crossing of arcs has been defined to match this geometrical picture in the obvious way.

\begin{defn} (\cite[Definition 3.2]{PJ11})
Let $\arcs$ be a set of arcs. If, for each integer $n$, there are only finitely many arcs in $\arcs$ which end in $n$, then $\arcs$ is said to be locally finite. A left (resp. right) fountain of $\arcs$ is an integer $n$ for which there are infinitely many arcs of the form $(m,n)$ (resp. $(n,m)$) in $\arcs$. A fountain of $\arcs$ is an integer $n$ which is both a left and a right fountain of $\arcs$. Finally, $\arcs$ is said to be non-crossing if $\arcs$ contains no pair of crossing arcs. 
\end{defn}



We can identify coordinate pairs $(m,n)$ such that $n-m \geq 2$ with either indecomposable objects of $\cal D$ or with arcs. Hence a set $\arcs$ of arcs gives a collection of indecomposable objects of $\cal D$, and add of these gives a subcategory $\cal A$ of $\cal D$. This sets up a bijection between sets of arcs and subcategories of $\cal D$.

Now let us describe the layout of this paper.

In Section~\ref{Definitionsarcs}, we give a few definitions. 





In Section~\ref{functoriallyfinitesubcategories}, we consider the precovering and preenveloping properties of a subcategory of $\cal D$ and show the following.

\begin{mythm}
Let $\cal A$ be a subcategory of $\cal D$ and let
$\arcs$ be the corresponding set of arcs.
Then $\cal A$ is precovering if and only if each right fountain of $\arcs$ is in fact a fountain.
\end{mythm}

\begin{mythm}
Let $\cal A$ be a subcategory of $\cal D$ and let
$\arcs$ be the corresponding set of arcs.
Then $\cal A$ is preenveloping if and only if each left fountain of $\arcs$ is in fact a fountain.
\end{mythm}

They will be proved in Theorem~\ref{precoveringfountain} and in Theorem~\ref{precoveringfountaindual} respectively.

In Section~\ref{configurationsofarcs}, we will need the following definitions. The first two describe some constraints on a set of arcs $\arcs$.

\begin{defn}   \label{config1}
A set of arcs $\arcs$ is said to satisfy \emph{condition (i)} if, for each pair of crossing arcs $(a, b)$ and $(c, d)$ in $\arcs$, those of the pairs $(a,c)$, $(c,b)$, $(b,d)$ and $(a,d)$ which are arcs belong to $\arcs$ (for instance, $(a,c)$ is only an arc if $c-a \geq 2$). 
\end{defn}

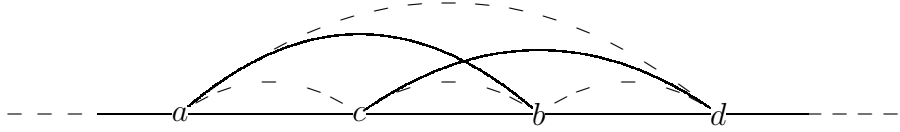
\[
  \xymatrix @-2.5pc @! {
    \rule{0ex}{6.5ex}\ar@{--}[r]
    &\ar@{-}[r]
&*{a} \ar@{-}[r]\ar@/^2.5pc/@{-}[rrrr]  \ar@/^1.0pc/@{--}[rr] \ar@/^3.5pc/@{--}[rrrrrr]
&*{}\ar@{-}[r]& *{c} \ar@{-}[rr]  \ar@/^1.0pc/@{--}[rr] \ar@/^2pc/@{-}[rrrr]
& \ar@{-}[r] 
& *{b}   \ar@{-}[rr] \ar@/^1.0pc/@{--}[rr]
& \ar@{-}[r] 
& *{d} \ar@{-}[r] 
& *{}\ar@{--}[r]
&*{}
                    }
\]

\begin{defn}  \label{config2}
A set of arcs  $\arcs$ is said to satisfy \emph{condition (ii)} if it has the following property: 
if $\textit{a}$ is a left fountain but not a right fountain of $\arcs$, $\textit{b}$ is a right fountain but not a left fountain of $\arcs$ and $b-a \geq 2$, then the arc $(a,b)$ is in $\arcs$.
\end{defn}

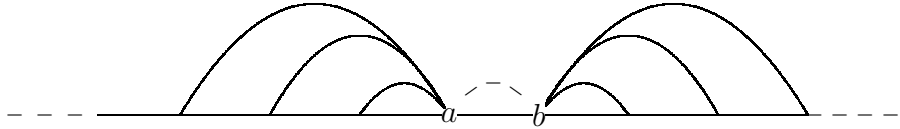
\[
  \xymatrix @-2.5pc @! {
    \rule{0ex}{6.5ex}\ar@{--}[r]
    &\ar@{-}[r]
&*{} \ar@{-}[r]  \ar@/^3.5pc/@{-}[rrr]    
&*{}\ar@{-}[r]& *{} \ar@{-}[r]   
& *{a} \ar@{-}[r]  \ar@/_2.5pc/@{-}[ll] \ar@/_1.0pc/@{-}[l] \ar@/^1.0pc/@{--}[r] 
& *{b}   \ar@{-}[rr]   \ar@/^3.5pc/@{-}[rrr]    \ar@/^2.5pc/@{-}[rr]  \ar@/^1.0pc/@{-}[r]  
& *{} \ar@{-}[r] 
& *{} \ar@{-}[r] 
& *{}\ar@{--}[r]
&*{}
                    }\]





Using these, we establish the following main result. 

\begin{mythm}
Let $\cal X$ be a subcategory of $\cal D$ and let $\arcsX$ be the corresponding set of arcs. Then the following conditions are equivalent.
\begin{itemize}
\item[(i)] $\arcsX$ satisfies condition (i) and condition (ii), and each right fountain of $\arcsX$ is in fact a fountain,
\item [(ii)] The subcategory $\cal X$ is precovering and is closed under extensions,
\item [(iii)] $(\cal X, \cal Y)$ is a torsion theory for some subcategory $\cal Y$ of $\cal D$. 
\end{itemize}
\end{mythm}

Note that in (iii), ${\cal Y} = \{ y$ $\in$ ${\cal D}  {\mid} ({\cal X}, y)=0 \}$. This will be proved in Theorem~\ref{torsiontheoryconfiguration1}.

Finally in Section~\ref{examplesofarcs}, we give a few examples, characterizing all t-structures and co-t-structures in $\cal D$.

\section{Definitions}\label{Definitionsarcs}
In this section we give a few definitions and state a few observations.

\begin{defn}
Let $\cal U$ be a subcategory of $\cal D$. Define ${\cal U}^{\perp} = \{ d$ $\in$ ${\cal D}  {\mid} (u,d)=0$ for all $u$ in $\cal U \}$ and $^{\perp}{\cal U} = \{ d$ $\in$ ${\cal D}  {\mid} (d,u)=0$  for all $u$ in $\cal U \}$.
\end{defn}

\begin{defn}
Let $\arcsU$ be a set of arcs. Define $\mathrm{ort}(\arcsU)$ to be the set of arcs $\mathrm{ort}(\arcsU)$ = $\{ \arcsd \mid \arcsd$ does not cross any arcs in $\arcsU \}$. 
\end{defn}








Now let $\cal X$ and $\cal Y$ be two subcategories of $\cal D$ such that $\cal X$ = $^{\perp}{\cal Y}$ and $\cal Y = {\cal X}^{\perp}$. Let $\arcsX$ and $\arcsY$ be the corresponding sets of arcs respectively. By \cite[Lemma 3.6]{PJ11}, we have $\arcsX = \{ \arcsd  \mid \Sigma \arcsd$ does not cross any arcs in $\arcsY \}$ and  $\arcsY = \{ \arcsd  \mid \Sigma^{-1} \arcsd$ does not cross any arcs in $\arcsX \}$.

Let $\Sigma^{-1} \arcsY = \arcsW$. We can now rewrite the above as 
$\arcsX = \{ \arcsd  \mid \arcsd$  does not cross any arcs in $\arcsW \} = \mathrm{ort}(\arcsW)$ and $\arcsW = \{ \arcsd  \mid \arcsd  $ does not cross any arcs in $\arcsX \} = \mathrm{ort} (\arcsX)$.


\begin{lem}\label{ort4}
Let ${\cal X}$ and $\arcsX$ be defined as above. Then ${\cal X}$ =  $^{\perp}({\cal X}^{\perp})$ is equivalent to $\arcsX = \mathrm{ort \, ort} (\arcsX)$.
\end{lem}

\begin{proof}
This is obvious by the above description. 
\end{proof}

\section{Precovering (preenveloping) subcategories}\label{functoriallyfinitesubcategories}
In this section we 
characterize precovering and preenveloping subcategories in terms of their corresponding sets of arcs. 

Let us recall the following property regarding the morphisms in the category $\cal D$. Let $x$ be an indecomposable object of $\cal D$. See \cite[Definition 2.1]{PJ11} for the definitions of the regions $\H^{-}(x)$ and $\H^{+}(x)$ in the AR quiver of $\cal D$, which are sketched as follows.

\[
\xymatrix @-4.5pc @! {
    &&&*{} &&&&&&&& *{}&& \\
    &&&& *{} \ar@{--}[ul] & & & & & & *{} \ar@{--}[ur] \\
    &*{}&& H^-(x) & & & & & & & & H^+(x) && *{}\\
    &&*{}\ar@{--}[ul]& & & & {\Sigma x \hspace{3ex}} \ar@{-}[ddll] \ar@{-}[uull] & {x} & {\hspace{3ex} \Sigma^{-1} x} \ar@{-}[ddrr] \ar@{-}[uurr]& & &&*{}\ar@{--}[ur]&\\ 
    && \\
    *{}\ar@{--}[r]&*{} \ar@{-}[rrr] && &  *{} \ar@{-}[uull] \ar@{-}[rrrrrr]& & & & & &  *{} \ar@{-}[uurr]\ar@{-}[rrr]&&&*{}\ar@{--}[r]&*{}\\
           }
\]

We write $\H(x) = \H^-(x) \cup \H^+(x)$.

\begin{property} (\cite[Corollary 2.3]{PJ11}) \label{morphismsinAinfty}
Let $x$ and $y$ be indecomposable objects of $\cal D$.  Then the following are
equivalent.
\begin{itemize}
\item[(i)] $(x,y) \neq 0$,
\item[(ii)] $(x,y) = k$,
\item[(iii)] $y \in \H(\Sigma x)$,
\item[(iv)] $x \in \H(\Sigma^{-1}y)$.
\end{itemize}
\end{property}

Now we are ready to give the two theorems in this section.

\begin{thm}\label{precoveringfountain}
Let $\cal A$ be a subcategory of $\cal D$ and let $\arcs$ be the corresponding set of arcs. Then $\cal A$ is precovering if and only if each right fountain of $\arcs$ is in fact a fountain.
\end{thm}

\begin{proof}

Suppose $\cal A$ is precovering. A right fountain of $\arcs$ is an integer $n$ for which there are infinitely many arcs of
the form $(n,p)$ in $\arcs$.  The corresponding collection $P$ of indecomposable objects of $\cal A$ lie on a diagonal half line $r$ in the AR quiver of $\cal D$.  The following sketch shows $r$ along with some of the indecomposable objects $a_i$ in $P$ (indicated by the black dots) and, in dotted lines, their respective regions $\H(\Sigma a_i)$.

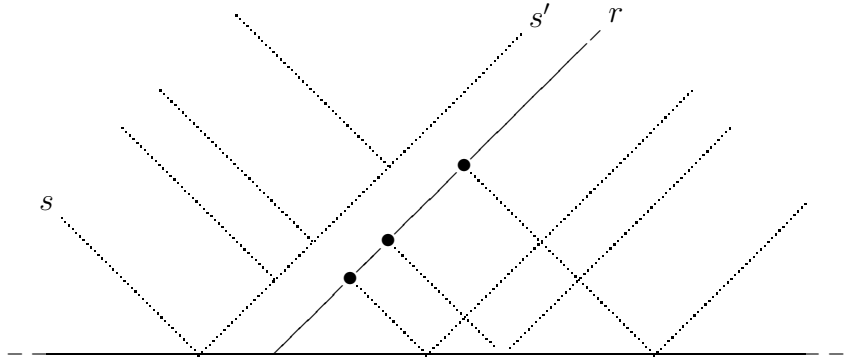
\[\xymatrix @-2.05pc @! {
    &&&&&&&*{}&&&&&&&&{s^{\prime}}&&{r} \\
    &&&&&&&&&&&&&&*{}&&*{}\ar@{--}[ur] \\
    &&&&&*{}&&&&&&&&*{}&&*{}&&&&*{}\\
    &&&&*{}&&&&&&&&&&*{}&&&&&*{}&*{}\\
    &&&&&*{} &&&&&*{}& *{} \ar@{.}[uuuullll]& *{} & *{\bullet} \ar@{-}[uuurrr]\ar@{.}[dddddrrrrr] &&&&*{}&&*{} \\
    &&{s}&& & & *{} & & & & & & *{} &&&&*{}&*{}&*{}&&&&*{}\\
    &&&&&  & & & & *{}\ar@{.}[uuuullll]& & *{\bullet} \ar@{-}[uurr] \ar@{.}[dddrrr]& & && *{}\\
    &&&  &  &  &    & & *{} \ar@{.}[uuuullll] & & *{\bullet} \ar@{.}[ddrr] \ar@{-}[ur]& & & &*{}&&&&&&*{}\\ 
    &&&& &  & &  & & &  &&&& \\
    &*{}\ar@{--}[r]&*{}  \ar@{-}[rrrr] && & & *{} \ar@{.}[uuuullll]\ar@{.}[uuuuuuuuurrrrrrrrr]\ar@{-}[rr]& *{} & *{} \ar@{-}[uurr]\ar@{-}[rrrr]& & & & *{} \ar@{.}[uuuuuuurrrrrrr]\ar@{-}[rrrrrrrrrr]&&\ar@{.}[uuuuuurrrrrr]&*{}&*{}&&*{}\ar@{.}[uuuurrrr]&&&&*{}\ar@{--}[r]&*{}\\
           }\]
We want to show that $n$ is also a left fountain, that is, there are infinitely many arcs in $\arcs$ of the form $(m,n)$.  This is the same as showing that there are infinitely many indecomposable objects of $\cal A$ which are on the half line $s$. Consider an object $x$ on the half line $s$ and its region $\H(\Sigma^{-1} x)$ in the following diagram. 

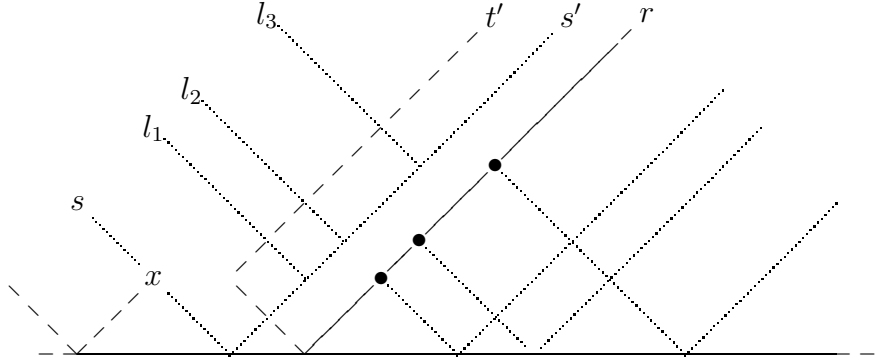
\[\xymatrix @-2.05pc @! {
    &&&&&&&*{l_3}&&&&&&{t^{\prime}}&&{s^{\prime}}&&{r} \\
    &&&&&&&&&&&&&&*{}&&*{}\ar@{--}[ur] \\
    &&&&&*{l_2}&&&&&&&&*{}&&*{}&&&&*{}\\
    &&&&*{l_1}&&&&&&&&&&*{}&&&&&*{}&*{}\\
    &&&&&*{} &&&&&*{}& *{} \ar@{.}[uuuullll]& *{} & *{\bullet} \ar@{-}[uuurrr]\ar@{.}[dddddrrrrr] &&&&*{}&&*{} \\
    &&{s}&& & & *{} & & & & & & *{} &&&&*{}&*{}&*{}&&&&*{}\\
    &&&&&  & & & & *{}\ar@{.}[uuuullll]& & *{\bullet} \ar@{-}[uurr] \ar@{.}[dddrrr]& & && *{}\\
    &&&  & x  \ar@{--}[ddll]  \ar@{.}[uull]&  &   \ar@{--}[uuuuuuurrrrrrr] \ar@{--}[ddrr] &  & *{} \ar@{.}[uuuullll] & & *{\bullet} \ar@{.}[ddrr] \ar@{-}[ur]& & & &*{}&&&&&&*{}\\ 
    &&&& &  & &  & & &  &&&& \\
    &*{}\ar@{--}[r]&*{} \ar@{--}[uull] \ar@{-}[rrrr] && & & *{} \ar@{.}[uull]\ar@{.}[uuuuuuuuurrrrrrrrr]\ar@{-}[rr]& *{} & *{} \ar@{-}[uurr]\ar@{-}[rrrr]& & & & *{} \ar@{.}[uuuuuuurrrrrrr]\ar@{-}[rrrrrrrrrr]&&\ar@{.}[uuuuuurrrrrr]&*{}&*{}&&*{}\ar@{.}[uuuurrrr]&&&&*{}\ar@{--}[r]&*{}\\
           }\]
           
Let $\beta: b \rightarrow x$ be an $\cal A$-precover. We write $b= b_1 \oplus \cdots \oplus b_q$ and $\beta = \begin{pmatrix} \beta_1,  \ldots,  \beta_q  \end{pmatrix}$. We can assume that the morphism $\beta: b \rightarrow x$ is non-zero on each direct summand $b_j$ of $b$. In particular, each direct summand $b_j$ of $b$ belongs to $\H(\Sigma^{-1}x)$. 

It is clear that $x$ is in $\H^-(\Sigma a_i)$ for all the $a$ in $P$. Pick a non-zero morphism $\alpha: a \rightarrow x$. The morphism factors through $\beta$ as shown in the following diagram.

\[\xymatrix @-3.5pc @! {  
a  \ar[drr]^{\alpha} \ar[d]_{\gamma} &  &   \\   
b_1 \oplus \cdots \oplus b_q  \ar[rr]_{\beta} & & x \\ 
 }\]

We write $\gamma = \begin{pmatrix} \gamma_{1} \\ \vdots  \\ \gamma_{q} \end{pmatrix}$. Since $\alpha$ is non-zero and $\alpha = \beta_1 \gamma_{1} + \cdots + \beta_q \gamma_{q}$, there is a term $\beta_k \gamma_{k}$ which is non-zero. Hence $\gamma_{k}: a \rightarrow b_k$ is non-zero so $b_k$ is in $\H(\Sigma a)$. 

Therefore $b_k$ can only be on the half line $s$ above $x$, or in the infinite region inside $\H^+(\Sigma^{-1} x)$ bounded by 
$t'$, $s'$ and $l_i$ for some $i$. 
But it cannot be the latter because then $\beta_k$ and $\gamma_{k}$ would both be backward morphisms whence $\beta_k \gamma_{k}$ would be zero. Therefore $b_k$ can only be on the half line $s$ above $x$. Repeating the argument, consider an object $x_1$ on the half line $s$ above $b_k$. We can then reveal another object $c$ in $\cal A$ on the half line $s$ above $x_1$. Since we can continue in this way indefinitely, there are infinitely many indecomposable objects of $\cal A$ which are on the half line $s$.


Now suppose each right fountain of $\arcs$ is in fact a fountain. 

In the following diagram, consider an object $x$ and its region $\H(\Sigma^{-1} x)$ in the AR quiver of $\cal D$. As indicated, we also introduce half lines (which start from the bottom line) of the form $r_i$ or $s_i$, where $i \in \mathbb{N}$ is a variable. Suppose the region $\H^{-}(\Sigma^{-1} x)$ (resp. $\H^{+}(\Sigma^{-1} x)$) has boundary half lines $s_0$ and $s_m$ (resp. $r_0$ and $r_m$). Then each line $s_i$ (resp. $r_i$), for $0 \le i \le m$, must pass through the region $\H^{-}(\Sigma^{-1} x)$ (resp. $\H^{+}(\Sigma^{-1} x)$) and is parallel to the boundary lines of the region $\H^{-}(\Sigma^{-1} x)$ (resp. $\H^{+}(\Sigma^{-1} x)$).

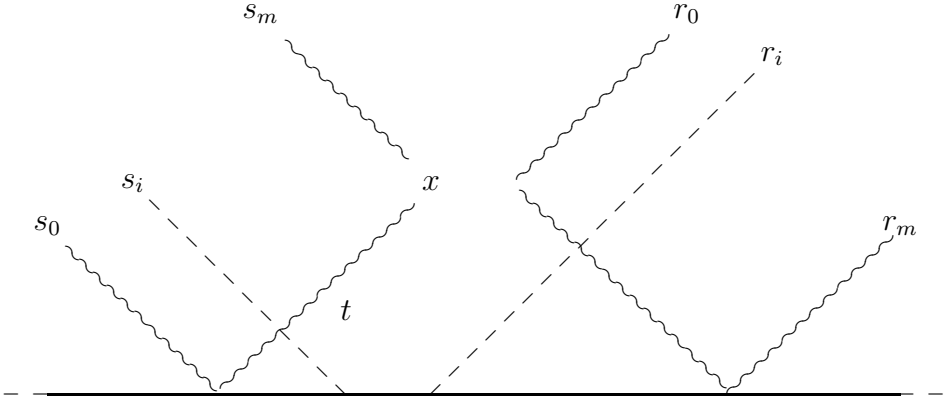
\[
\vcenter{
  \xymatrix @-2.3pc @! {
    &&&&&{}&s_m&&&&*{}&&&&&&r_0&&&& \\
    &&&&&&&&&&&&&&&&&&{r_i}&& \\
    &&&&&&&&&*{}&&&&&&&&*{} &&{}&&&&*{}\\
    &&&&{}&&&&*{}&&&&&&&&&&*{}&&&&\\
    &&&{s_i}&&  &&&&& {x} \ar@{~}[dddddlllll]\ar@{~}[uuuullll]&&*{}\ar@{~}[dddddrrrrr]\ar@{~}[uuuurrrr]&& *{}  & *{}  *{} & *{} &  &&&&*{}& \\
    &s_0&&&*{}  & &  && &  & & &&& & & *{} &&&&*{}&*{r_m}&*{}\\
    &&&& &*{} &&*{}&& & & & & *{}& &   & & && *{}&&&\\
    &&&&&&*{}  &&*{t}& & & & *{}  & & *{} & & & &*{}&&&&*{}\\ 
    &&&& & & &  &&  &&& & & &&&& \\
    *{}\ar@{--}[r]&*{}\ar@{-}[rrrr]&&&&*{} \ar@{~}[uuuullll]\ar@{-}[rrrrr] &&& *{} \ar@{--}[uuuuulllll] &  & *{} \ar@{--}[uuuuuuuurrrrrrrr]\ar@{-}[rr]& *{} & *{} \ar@{-}[rrrr]& & & &  *{} \ar@{-}[rrrrr]&*{}\ar@{~}[uuuurrrr]&&*{}&&*{}\ar@{--}[r]&*{}\\
           }
}
\]

Let $S$ be the intersection of $\H^{-}(\Sigma^{-1} x)$ and the objects of $\cal A$. On each line $s_i$, $0 \le i \le m$, consider the first object $a_i$ in $\cal A$ which lies above the line segment $t$. Denote by $a_s$ the direct sum of all the $a_i$ and consider the canonical morphism $a_s \rightarrow x$. By \cite[Lemma 2.5]{PJ11}, each morphism $a \rightarrow x$ with $a$ in $S$ factors as $a \rightarrow a_s \rightarrow x$. 

For example, in the following diagram, we have $m = 5$ and the 
circles and bullets indicate the first few objects of $\cal A$ on each line $s_i$. The object $a_s$ described above is the direct sum of the objects indicated by circles $\circ$.

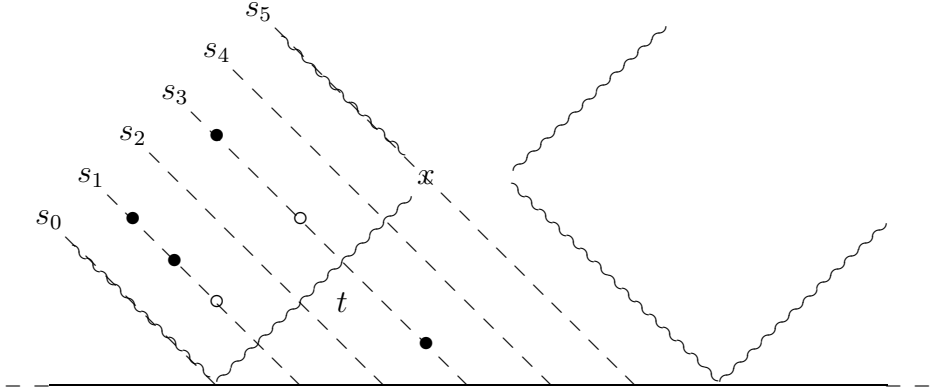
\[
\vcenter{
  \xymatrix @-2.05pc @! {
    &&&&&&&{s_5}&&&&&*{}&&&&&&&&&& \\
    &&&&&&s_4&&&&&&&&&&&&&{}&& \\
    &&&&&s_3&&&&&*{}&&&&&&&&*{} &&{}&&&&*{}\\
    &&&&s_2&&*{\bullet}&&&&&&&&&&&&&*{}&&&&\\
    &&&s_1&&&  &&&&& {x} \ar@{~}[dddddlllll]\ar@{~}[uuuullll]&&*{}\ar@{~}[dddddrrrrr]\ar@{~}[uuuurrrr]&& *{}  &   *{} & *{} &  &&&&*{}& \\
    &&s_0&&*{\bullet}&  & &  &*{\circ}& &  & & &&& & & *{} &&&&*{}&*{}&*{}\\
    &&&&& *{\bullet}& &&*{}&& & & & & *{}& &   & & && *{}&&&\\
    &&&&&&*{\circ}&  &&*{t}& & & & *{}  & & *{} & & & &*{}&&&&*{}\\ 
    &&&&& & & &  &&  &*{\bullet}&& & & &&&& \\
    &*{}\ar@{--}[r]&*{}\ar@{-}[rrrr]&&&&*{} \ar@{~}[uuuullll] \ar@{--}[uuuullll] \ar@{-}[rrrrr] & *{} & *{}\ar@{--}[uuuuulllll]& *{}  &  *{}\ar@{--}[uuuuuullllll] & *{} \ar@{-}[rr]  & *{} *{} \ar@{--}[uuuuuuulllllll] & *{} \ar@{-}[rrrr]  & *{}\ar@{--}[uuuuuuuullllllll] &  & *{}\ar@{--}[uuuuuuuuulllllllll] &  *{} \ar@{-}[rrrrr]&*{}\ar@{~}[uuuurrrr]&&*{}&&*{}\ar@{--}[r]&*{}\\
           }
}
\]


Let $R$ be the intersection of $\H^{+}(\Sigma^{-1} x)$ and the objects of $\cal A$. We want to find an object $a_r$ in $\cal A$ with a non-zero morphism $a_r \rightarrow x$ such that each morphism $a \rightarrow x$ with $a$ in $R$ factors as $a \rightarrow a_r \rightarrow x$.

Suppose $R$ is finite. Then 
we can let $a_r$ be the direct sum of the objects in $R$. Otherwise $R$ is infinite. Since there are only finitely many $r_i$, therefore there is a line $r_j$ which contains infinitely many objects in $\cal A$ with a non-zero morphism to $x$. Let $J \subseteq \{0, \ldots, m\}$ consist of the $j$ such that the line $r_j$ contains infinitely many objects in $\cal A$ with a non-zero morphism to $x$. 

Now for each $j$ in $J$, the line $r_j$ corresponds to a right fountain of $\arcs$, and by assumption, each right fountain is also a fountain. Hence the corresponding line $s_j$ contains infinitely many objects in $\cal A$ with a non-zero morphism to $x$. Among the $s_j$ with $j$ in $J$, take the line $s_q$ which is closest to the boundary line $s_m$, and then consider any object which is in both $\cal A$ and $\H^{-}(\Sigma^{-1} x)$ on that line. By \cite[Lemma 2.7]{PJ11}, this object plays the role of $a_r$. 

For example, in the following diagram, $m=5$ and $J = \{1, 2\}$. The line $s_q$ described above is the line $s_2$ here. The 
bullets on the $r_i$ indicate the first few objects of $\cal A$ in $\H^{+}(\Sigma^{-1} x)$. The 
bullets and the circle on the $s_i$ indicate the first few objects of $\cal A$ on those half lines.  The object $a_r$ described above is indicated by the circle $\circ$ (one of the infinitely many choices).

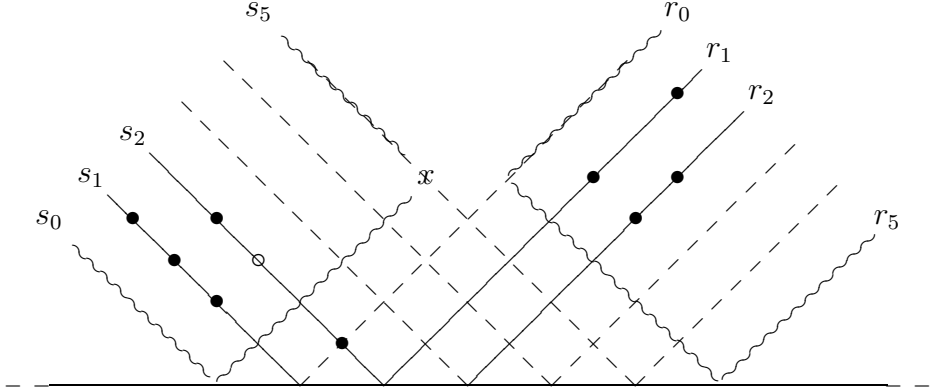
\[
\vcenter{
  \xymatrix @-2.05pc @! {
    &&&&&&&{s_5}&&&&&*{}&&&&&r_0&&&&& \\
    &&&&&&&&&&&&&&&&&&r_1&{}&& \\
    &&&&&&&&&&*{}&&&&&&&*{\bullet}& &r_2&{}&&&&*{}\\
    &&&&s_2&{}&&&&*{}&&&&&&&&&*{}&*{}&&&&\\
    &&&s_1&{}&*{}& *{} &*{}&&&& {x} \ar@{~}[dddddlllll]\ar@{~}[uuuullll]&&*{}\ar@{~}[dddddrrrrr]\ar@{~}[uuuurrrr]&& *{\bullet}  &   *{} & *{\bullet}&   &&&&*{}& \\
    &&s_0&&*{\bullet}&*{}  & *{\bullet} &  &*{}& &  & & &&& &*{\bullet} &  &&&&*{}&*{r_5}&*{}\\
    &&&&& *{\bullet}&*{} &*{\circ}&&& & & & & *{}& &   & & && *{}&&&\\
    &&&&&&*{\bullet}&*{}  &&& & & & *{}  & & *{} & & & &*{}&&&&*{}\\ 
    &&&&& & & &  &*{\bullet}&  &&& & & &&&& \\
    &*{}\ar@{--}[r]&*{}\ar@{-}[rrrr]&&&  &*{} \ar@{~}[uuuullll]\ar@{-}[rrrrr] &&*{}\ar@{-}[uuuuulllll] \ar@{--}[uuuuuuuurrrrrrrr]& *{}  &  *{}\ar@{-}[uuuuuullllll] \ar@{-}[uuuuuuuurrrrrrrr]  & *{} \ar@{-}[rr]& *{}  *{}\ar@{--}[uuuuuuulllllll] \ar@{-}[uuuuuuurrrrrrr] & *{} \ar@{-}[rrrr]&  *{}\ar@{--}[uuuuuuuullllllll] \ar@{--}[uuuuuurrrrrr]& & *{}\ar@{--}[uuuuuuuullllllll] \ar@{--}[uuuuurrrrr]&  *{} \ar@{-}[rrrrr]&*{}\ar@{~}[uuuurrrr]&&*{}&&*{}\ar@{--}[r]&*{}\\
           }
}
\]

Finally, an $\cal A$-precover of $x$ can be obtained as $a_r \oplus a_s \rightarrow x$. 

\end{proof}

We also have the dual of the above theorem.

\begin{thm}\label{precoveringfountaindual}
Let $\cal A$ be a subcategory of $\cal D$ and let $\arcs$ be the corresponding set of arcs. Then $\cal A$ is preenveloping if and only if each left fountain of $\arcs$ is in fact a fountain.
\end{thm}

\begin{proof}
Similar.
\end{proof}

\section{Torsion theories}  \label{configurationsofarcs}
In this section, let $\arcs$ be a set of arcs. We first have a sequence of results regarding $\arcs$, and then we will give 
a checkable condition equivalent to $\arcsX = \mathrm{ort \, ort} (\arcsX)$, see Lemma~\ref{configurationequivelance}. 
Finally, we will give the main theorem of this paper, see Theorem~\ref{torsiontheoryconfiguration1}.



\begin{lem}  \label{configurationimplication}
Suppose $\mathrm{ort \, ort} \, {\arcs} = \arcs$. Then $\arcs$ satisfies conditions (i) and (ii).
\end{lem}

\begin{proof}
To see that $\arcs$ satisfies condition (i), consider the diagram in Definition~\ref{config1}. The arc $(a,c)$ is in $\mathrm{ort \, ort} \, {\arcs}$, otherwise there would have to be an arc $(m,n)$ in $\mathrm{ort} \, \arcs$ which crossed $(a,c)$, but the arc $(m,n)$ cannot be in $\mathrm{ort} \, \arcs$ since it crosses either the arc $(a,b)$ or the arc $(c,d)$ in $\arcs$. Therefore the arc $(a,c)$ is in $\mathrm{ort \, ort} \, {\arcs} = \arcs$. The rest is similar. 

To see that $\arcs$ satisfies condition (ii), consider the diagram in Definition~\ref{config2}. Suppose the arc $(a,b)$ is not in  $\mathrm{ort \, ort} \, {\arcs}$. Then there is an arc $(m,n)$ in $\mathrm{ort} \, \arcs$ which crosses $(a,b)$, that is, $m < a < n < b$ or $a < m < b < n$. For the first case the arc $(m,n)$ cannot however be in $\mathrm{ort} \, \arcs$, since there is always an arc $(q,a)$ with $q < m$ in $\arcs$ which crosses $(m,n)$, and similarly for the second case. Therefore the arc $(a,b)$ is in $\mathrm{ort \, ort} \, {\arcs} = \arcs$.
\end{proof}

\begin{lem} \label{configurationleftright}
Suppose $\arcs$ satisfies condition (i) and suppose there are only finitely many (but not zero) arcs in $\arcs$ that end in $a$. Suppose there are both arcs going to the left and arcs going to the right from $a$. If $(p,a)$ is the longest arc in $\arcs$ going to the left from $a$ and $(a,q)$ is the longest arc in $\arcs$ going to the right from $a$, then $(p,q)$ is an arc in $\mathrm{ort} \, \arcs$.
\end{lem}
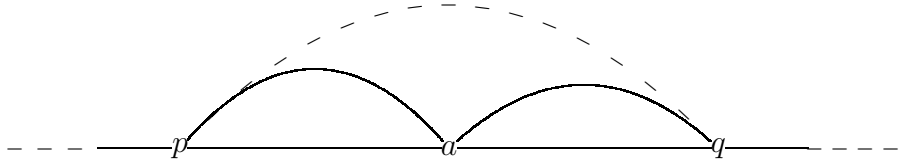
\[
  \xymatrix @-2.5pc @! {
   \\
    \rule{0ex}{6.5ex}\ar@{--}[r]
    &\ar@{-}[r]
&*{p} \ar@{-}[r]  \ar@/^2.5pc/@{-}[rrr]    \ar@/^4.5pc/@{--}[rrrrrr]
&*{}\ar@{-}[r]& *{} \ar@{-}[r]   
& *{a} \ar@{-}[r]  *{} \ar@/^2pc/@{-}[rrr]
& *{}   \ar@{-}[rr] 
& *{} \ar@{-}[r] 
& *{q} \ar@{-}[r] 
& *{}\ar@{--}[r]
&*{}
                    }
\]
\begin{proof}
There are no arcs $(m,n)$ in $\arcs$ with $m < p$, $p < n  < a$, otherwise $(m,a)$ would be in $\arcs$ by condition (i), contradicting that $(p,a)$ is the longest arc in $\arcs$ going to the left from $a$. There are also no arcs $(m,n)$ in $\arcs$ with $p < m < a$, $q < n $, otherwise $(a,n)$ would be in $\arcs$ by condition (i), contradicting that $(a,q)$ is the longest arc in $\arcs$ going to the right from $a$. 

Similarly, there are no arcs $(m,n)$ in $\arcs$ with $a < m < q$, $q < n $, and there are also no arcs $(m,n)$ in $\arcs$ with $m < p$, $a < n < q $. By construction there are no arcs $(m,a)$ in $\arcs$ with $m < p$, and there are also no arcs $(a,n)$ in $\arcs$ with $q < n$. 

Combining all these shows that there are no arcs in $\arcs$ crossing $(p,q)$ so $(p,q)$ has to be in $\mathrm{ort} \, \arcs$.
\end{proof}

\begin{lem}   \label{configurationleft}
Suppose $\arcs$ satisfies condition (i) and suppose there are only finitely many (but not zero) arcs in $\arcs$ that end in $a$. Suppose there are only arcs going to the left from $a$. If $(p,a)$ is the longest arc in $\arcs$ going to the left from $a$, then $(p,a+1)$ is in $\mathrm{ort} \, \arcs$.
\end{lem}
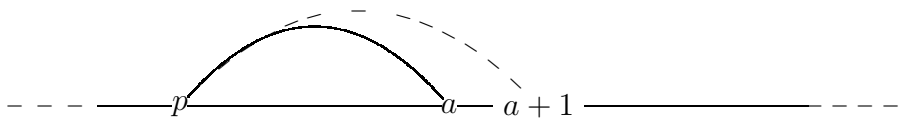
\[
  \xymatrix @-2.5pc @! {
    \rule{0ex}{6.5ex}\ar@{--}[r]
    &\ar@{-}[r]
&*{p} \ar@{-}[r]  \ar@/^2.5pc/@{-}[rrr]    \ar@/^3.0pc/@{--}[rrrr]
&*{}\ar@{-}[r]& *{} \ar@{-}[r]   
& *{a} \ar@{-}[r]  
& *{\ a+1 \ }   \ar@{-}[rr] 
& *{} \ar@{-}[r] 
& *{} \ar@{-}[r] 
& *{}\ar@{--}[r]
&*{}
                    }
\]

\begin{proof}
There are no arcs $(m,n)$ in $\arcs$ with $m < p$, $p < n  < a$, otherwise $(m,a)$ would be in $\arcs$ by condition (i), contradicting that $(p,a)$ is the longest arc in $\arcs$ going to the left from $a$. There are also no arcs $(m,n)$ in $\arcs$ with $p < m  < a$, $n > a+1$, otherwise $(a,n)$ would be in $\arcs$ by condition (i), contradicting that there are no arcs in $\arcs$ going to the right from $a$. By construction there are no arcs $(m,a)$ in $\arcs$ with $m < p$ and it is a condition that there are no arcs $(a,n)$ in $\arcs$. Therefore $(p,a+1)$ has to be in $\mathrm{ort} \, \arcs$.

\end{proof}
     
\begin{lem}    \label{configurationright}
Suppose $\arcs$ satisfies condition (i) and suppose there are only finitely many (but not zero) arcs in $\arcs$ that end in $a$. Suppose there are only arcs going to the right from $a$. If $(a,p)$ is the longest arc in $\arcs$ going to the right from $a$, then $(a-1,p)$ is in $\mathrm{ort} \, \arcs$.  
\end{lem}
\[
  \xymatrix @-2.5pc @! {
    \rule{0ex}{6.5ex}\ar@{--}[r]
    &\ar@{-}[r]
&*{\ a-1 \ } \ar@{-}[r]       \ar@/^3.5pc/@{--}[rrrr]
&*{a}\ar@{-}[r] \ar@/^2.5pc/@{-}[rrr] & *{} \ar@{-}[r]  
& *{} \ar@{-}[r]  
& *{p}   \ar@{-}[rr] 
& *{} \ar@{-}[r] 
& *{} \ar@{-}[r] 
& *{}\ar@{--}[r]
&*{}
                    }
\]
     
\begin{proof}
Similar to Lemma~\ref{configurationleft}.
\end{proof}

\begin{remark} 
Suppose $(a,b)$ is an arc in $\mathrm{ort}^{2}(\arcs)$. Then there must be some arc in $\arcs$ which ends in $a$ (resp. $b$). Otherwise the arc $(a-1, a+1)$ (resp. $(b-1,b+1)$) is in $\mathrm{ort} \, \arcs$, which is a contradiction since it crosses $(a,b)$ in $\mathrm{ort}^{2}(\arcs)$.
\end{remark}
  
\begin{corol}   \label{corolconfigurationleft1}
Suppose $\arcs$ satisfies condition (i). Suppose that $(a,b)$ is an arc in $\mathrm{ort}^{2}(\arcs)$ and that there are only finitely many arcs in $\arcs$ that end in $a$. Then it is not possible that all arcs in $\arcs$ that end in $a$ are of the form $(m,a)$. 
\end{corol}

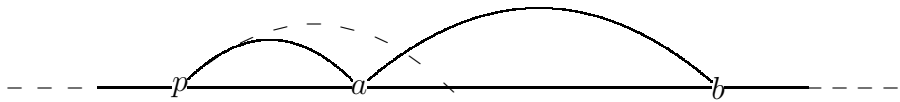
\[
  \xymatrix @-2.5pc @! {
    \rule{0ex}{6.5ex}\ar@{--}[r]
    &\ar@{-}[r]
&*{p} \ar@{-}[r]       \ar@/^1.5pc/@{-}[rr]  \ar@/^2.0pc/@{--}[rrr]
&*{}\ar@{-}[r]  & *{a} \ar@{-}[rrrr]   \ar@/^2.5pc/@{-}[rrrr]
& *{} 
& *{}   
& *{} 
& *{b} \ar@{-}[r] 
& *{}\ar@{--}[r]
&*{}
                    }
\]

\begin{proof}
Suppose all arcs in $\arcs$ that end in $a$ were of the form $(m,a)$ and let $(p,a)$ be the longest one. Then by Lemma~\ref{configurationleft} the arc $(p,a+1)$ is in ort $\arcs$, but this is a contradiction since it crosses $(a,b)$ in $\mathrm{ort}^{2}(\arcs)$. 
\end{proof}

\begin{corol}   \label{corolconfigurationleft2}
Suppose $\arcs$ satisfies condition (i). Let $(a,b)$ be an arc in $\mathrm{ort}^{2}(\arcs)$ and suppose there are only finitely many arcs in $\arcs$ that end in $a$. Then there is an arc $(a,m)$ in $\arcs$ with $b \leq m$.
\end{corol}

\begin{proof}
By Corollary~\ref{corolconfigurationleft1}, it is not possible that all arcs in $\arcs$ that end in $a$ are of the form $(m,a)$. Therefore we only have the remaining two cases. 
\begin{itemize}
  \item[(i)] Suppose there are arcs in $\arcs$ going to the left and going to the right from $a$. Let $(a,q)$ be the longest arc in $\arcs$ going to the right from $a$ and $(p,a)$ be the longest arc in $\arcs$ going to the left from $a$. By Lemma~\ref{configurationleftright}, $(p,q)$ is in $\mathrm{ort} \, \arcs$, so $q < b$ is not possible since then $(p,q)$ would cross $(a,b)$ in $\mathrm{ort}^{2}(\arcs)$. 
  
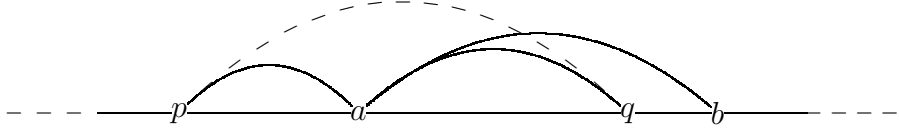
  \[
  \xymatrix @-2.5pc @! {
    \rule{0ex}{6.5ex}\ar@{--}[r]
    &\ar@{-}[r]
&*{p} \ar@{-}[r]       \ar@/^1.5pc/@{-}[rr]  \ar@/^3.5pc/@{--}[rrrrr]
&*{}\ar@{-}[r]  & *{a} \ar@{-}[r]   \ar@/^2.5pc/@{-}[rrrr] \ar@/^2.0pc/@{-}[rrr]
& *{} \ar@{-}[r] 
& *{}   \ar@{-}[r] 
& *{q} \ar@{-}[r] 
& *{b} \ar@{-}[r] 
& *{}\ar@{--}[r]
&*{}
                    }
\]

Therefore we have $b \leq q$. 

  \item[(ii)] Suppose there are only arcs in $\arcs$ going to the right from $a$. Let $(a,p)$ be the longest arc in $\arcs$ going to the right from $a$. By Lemma~\ref{configurationright}, $(a-1,p)$ is in $\mathrm{ort} \, \arcs$, so $p < b$ is not possible since then $(a-1,p)$ would cross $(a,b)$ in $\mathrm{ort}^{2}(\arcs)$. 
  
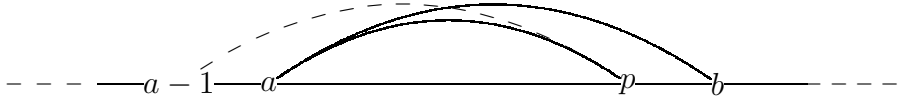
 \[
  \xymatrix @-2.5pc @! {
    \rule{0ex}{6.5ex}\ar@{--}[r]
    &\ar@{-}[r]
&*{a-1} \ar@{-}[r]       \ar@/^2.5pc/@{--}[rrrrr]
&*{a}\ar@{-}[r] \ar@/^2.0pc/@{-}[rrrr] \ar@/^2.5pc/@{-}[rrrrr]  & *{} \ar@{-}[r]  
& *{} \ar@{-}[r]  
& *{}   \ar@{-}[r] 
& *{p} \ar@{-}[r] 
& *{b} \ar@{-}[r] 
& *{}\ar@{--}[r]
&*{}
                    }
\]
     
Therefore we have $b \leq p$. 
  \end{itemize}
  
\end{proof}

\begin{corol}   \label{corolconfigurationright1}
Suppose $\arcs$ satisfies condition (i). Suppose that $(a,b)$ is an arc in $\mathrm{ort}^{2}(\arcs)$ and that there are only finitely many arcs in $\arcs$ that end in $b$. Then it is not possible that all arcs in $\arcs$ that end in $b$ are of the form $(b,m)$.
\end{corol}

\begin{proof}
Similar to Corollary~\ref{corolconfigurationleft1}.
\end{proof}

\begin{corol}   \label{corolconfigurationright2}
Suppose $\arcs$ satisfies condition (i). Let $(a,b)$ be an arc in $\mathrm{ort}^{2}(\arcs)$ and suppose there are only finitely many arcs in $\arcs$ that end in $b$. Then there is an arc $(m,b)$ in $\arcs$ with $m \leq a$.
\end{corol}

\begin{proof}
Similar to Corollary~\ref{corolconfigurationleft2}.
\end{proof}

\begin{lem}  \label{locallyfinite}
Suppose $\arcs$ satisfies condition (i). Let $(a,b)$ be an arc in $\mathrm{ort}^{2}(\arcs)$ and suppose that each of $a$ and $b$ is only an end point of finitely many arcs in $\arcs$. Then $(a,b)$ is in $\arcs$.
\end{lem}

\begin{proof}
By Corollary~\ref{corolconfigurationleft2}, there is an arc $(a,q)$ in $\arcs$ with $b \leq q$. On the other hand by Corollary~\ref{corolconfigurationright2}, there is an arc $(p,b)$ in $\arcs$ with $p \leq a$. If $q = b$  or $p=a$ already, then we are done. Assume otherwise that $q > b $ and $p < a$. But then $(a,b)$ is in $\arcs$ by condition (i).

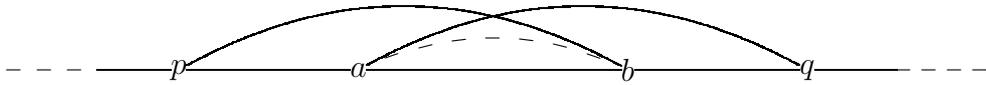
\[
  \xymatrix @-2.5pc @! {
    \rule{0ex}{6.5ex}\ar@{--}[r]
    &\ar@{-}[r]
&*{p} \ar@{-}[r]          \ar@/^2.0pc/@{-}[rrrrr] 
&*{}\ar@{-}[r]  & *{a} \ar@{-}[r] \ar@/^1.0pc/@{--}[rrr] \ar@/^2.0pc/@{-}[rrrrr] 
& *{} \ar@{-}[r] 
& *{}   \ar@{-}[r]  
& *{b} \ar@{-}[r] 
& *{} \ar@{-}[r] 
& *{q} \ar@{-}[r]
& *{}\ar@{--}[r]
&*{}
                    }
\]

\end{proof}

\begin{lem}  \label{fountain1}
Suppose $\arcs$ satisfies condition (i).
Let both $a$ and $b$ be right (resp. left) fountains of $\arcs$ with $b - a \geq 2$. Then $(a,b)$ is in $\arcs$.
\end{lem}

\begin{proof}
Suppose both $a$ and $b$ are right fountains of $\arcs$.
Choose an arc $(a,p)$ in $\arcs$ with $b < p$ and then choose an arc $(b,q)$ in $\arcs$ with $p < q$. Then $(a,b)$ is in $\arcs$ by condition (i). 
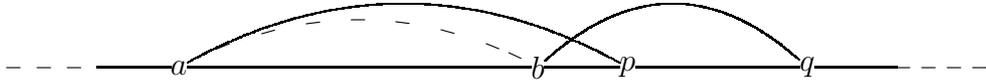
\[
  \xymatrix @-2.5pc @! {
    \rule{0ex}{6.5ex}\ar@{--}[r]
    &\ar@{-}[r]
&*{a} \ar@{-}[r]          \ar@/^2.0pc/@{-}[rrrrr] \ar@/^1.5pc/@{--}[rrrr] 
&*{}\ar@{-}[r]  & *{} \ar@{-}[r] 
& *{} \ar@{-}[r]  
& *{b}   \ar@{-}[r]  \ar@/^2.0pc/@{-}[rrr] 
& *{p} \ar@{-}[r] 
& *{} \ar@{-}[r] 
& *{q} \ar@{-}[r]
& *{}\ar@{--}[r]
&*{}
                    }
\]
The other case is similar.
\end{proof}

\begin{lem}  \label{fountain2}
Suppose $\arcs$ satisfies condition (i). Let $a$ be a right fountain of $\arcs$, $b$ be a left fountain of $\arcs$ with $b-a \geq 2$. Then $(a,b)$ is in $\arcs$.
\end{lem}

\begin{proof}
Choose an arc $(a,q)$ in $\arcs$ with $b < q$ and then choose an arc $(p,b)$ in $\arcs$ with $p < a$. Then $(a,b)$ is in $\arcs$ by condition (i).
\[
  \xymatrix @-2.5pc @! {
    \rule{0ex}{6.5ex}\ar@{--}[r]
    &\ar@{-}[r]
&*{p} \ar@{-}[r]          \ar@/^2.0pc/@{-}[rrrrr] 
&*{}\ar@{-}[r]  & *{a} \ar@{-}[r] \ar@/^1.0pc/@{--}[rrr] \ar@/^2.0pc/@{-}[rrrrr] 
& *{} \ar@{-}[r] 
& *{}   \ar@{-}[r]  
& *{b} \ar@{-}[r] 
& *{} \ar@{-}[r] 
& *{q} \ar@{-}[r]
& *{}\ar@{--}[r]
&*{}
                    }
\]

\end{proof}

\begin{lem}  \label{fountain3}
Suppose $\arcs$ satisfies condition (i). Let $(a,b)$ be an arc in $\mathrm{ort}^{2}(\arcs)$. Suppose $b$ is a left fountain of $\arcs$ and suppose there are only finitely many arcs in $\arcs$ that end in $a$. Then $(a,b)$ is in $\arcs$.
\end{lem}

\begin{proof}

By Corollary~\ref{corolconfigurationleft2}, there is an arc $(a,m)$ in $\arcs$ with $b \leq m$. If $m=b$ already then we are done. Otherwise choose an arc $(p,b)$ in $\arcs$ with $p<a$. Then $(a,b)$ is in $\arcs$ by condition (i).

\[
  \xymatrix @-2.5pc @! {
    \rule{0ex}{6.5ex}\ar@{--}[r]
    &\ar@{-}[r]
&*{p} \ar@{-}[r]          \ar@/^2.0pc/@{-}[rrrrr] 
&*{}\ar@{-}[r]  & *{a} \ar@{-}[r] \ar@/^1.0pc/@{--}[rrr] \ar@/^2.0pc/@{-}[rrrrr] 
& *{} \ar@{-}[r] 
& *{}   \ar@{-}[r]  
& *{b} \ar@{-}[r] 
& *{} \ar@{-}[r] 
& *{m} \ar@{-}[r]
& *{}\ar@{--}[r]
&*{}
                   }
\]

\end{proof}

\begin{lem}  \label{fountain4}
Suppose $\arcs$ satisfies condition (i). Let $(a,b)$ be an arc in $\mathrm{ort}^{2}(\arcs)$. Suppose $b$ is a right fountain of $\arcs$ and suppose there are only finitely many arcs in $\arcs$ that end in $a$. Then $(a,b)$ is in $\arcs$.
\end{lem}

\begin{proof}
By Corollary~\ref{corolconfigurationleft2}, there is an arc $(a,m)$ in $\arcs$ with $b \leq m$. If $m=b$ already then we are done. Otherwise choose an arc $(b,q)$ in $\arcs$ with $m<q$. Then $(a,b)$ is in $\arcs$ by condition (i).

\[
  \xymatrix @-2.5pc @! {
    \rule{0ex}{6.5ex}\ar@{--}[r]
    &\ar@{-}[r]
&*{a} \ar@{-}[r]          \ar@/^2.0pc/@{-}[rrrr] \ar@/^1.0pc/@{--}[rrr]
&*{}\ar@{-}[r]  & *{} \ar@{-}[r]  
& *{b} \ar@{-}[r]  \ar@/^2.0pc/@{-}[rrrr]
& *{m}   \ar@{-}[rr]  
& *{} \ar@{-}[r] 
& *{} \ar@{-}[r] 
& *{q} \ar@{-}[r]
& *{}\ar@{--}[r]
&*{}
                    }
\]

\end{proof}

\begin{lem}  \label{fountain5}
Suppose $\arcs$ satisfies condition (i). Let $(a,b)$ be an arc in $\mathrm{ort}^{2}(\arcs)$. Suppose $a$ is a right fountain of $\arcs$ and suppose there are only finitely many arcs in $\arcs$ that end in $b$. Then $(a,b)$ is in $\arcs$.
\end{lem}

\begin{proof}


Similar to Lemma~\ref{fountain3}.
\end{proof}

\begin{lem}  \label{fountain6}
Suppose $\arcs$ satisfies condition (i). Let $(a,b)$ be an arc in $\mathrm{ort}^{2}(\arcs)$. Suppose $a$ is a left fountain of $\arcs$ and suppose there are only finitely many arcs in $\arcs$ that end in $b$. Then $(a,b)$ is in $\arcs$.
\end{lem}

\begin{proof}


Similar to Lemma~\ref{fountain4}.
\end{proof}

Finally, we have the following lemma which is a recollection of the above lemmas.
\begin{lem}  \label{configurationequivelance}
(c.f. Lemma~\ref{ort4})
$\mathrm{ort \, ort} \, \arcs = \arcs$ if and only if $\arcs$ satisfies conditions (i) and (ii).
\end{lem}

\begin{proof}

(\emph{only if}) This is Lemma~\ref{configurationimplication}.
(\emph{if}) It is clear that $\arcs \subseteq \mathrm{ort \, ort} \,  \arcs$. Let $(a,b)$ be an arc in $\mathrm{ort}^{2}(\arcs)$. Suppose that each of $a$ and $b$ is only an end point of finitely many arcs in $\arcs$. Then $(a,b)$ is in $\arcs$ by Lemma~\ref{locallyfinite}. Otherwise suppose that both $a$ and $b$ are end points of infinitely many arcs in $\arcs$. Then $(a,b)$ is in $\arcs$ by Lemma~\ref{fountain1}, Lemma~\ref{fountain2} and condition (ii). Finally, suppose that precisely one of $a$ and $b$ is an end point of finitely many arcs in $\arcs$. Then $(a,b)$ is in $\arcs$ by Lemma~\ref{fountain3}, Lemma~\ref{fountain4}, Lemma~\ref{fountain5} and Lemma~\ref{fountain6}.
\end{proof}

Now we are ready to give the main theorem of this paper.


\begin{thm}\label{torsiontheoryconfiguration1}
Let $\cal X$ be a subcategory of $\cal D$ and let $\arcsX$ be the corresponding set of arcs. Then the following conditions are equivalent.
\begin{itemize}
\item[(i)] $\arcsX$ satisfies conditions (i) and (ii), and each right fountain of $\arcsX$ is in fact a fountain,
\item [(ii)] The subcategory $\cal X$ is precovering and is closed under extensions,
\item [(iii)] $(\cal X, \cal Y)$ is a torsion theory for some subcategory $\cal Y$ of $\cal D$. In particular, $\cal Y = {\cal X}^{\perp}$.
\end{itemize}
\end{thm}

\begin{proof}
(i) $\Rightarrow$  (ii): 
We see that $\arcsX$ satisfying conditions (i) and (ii) implies $\cal X$ being closed under extensions by Lemma~\ref{ort4} and Lemma~\ref{configurationequivelance}. Finally $\cal X$ is precovering if and only if each right fountain of $\arcsX$ is in fact a fountain by Theorem~\ref{precoveringfountain}.
(ii) $\Leftrightarrow$  (iii): This is true by \cite[Proposition 2.3]{OIYY 1}.
(iii) $\Rightarrow$  (i): Since $(\cal X, \cal Y)$ is a torsion theory, therefore ${\cal X}$ =  $^{\perp}({\cal X}^{\perp})$ and $\cal X$ is precovering. Therefore $\arcsX$ satisfies conditions (i) and (ii) by Lemma~\ref{ort4} and Lemma~\ref{configurationequivelance} and each right fountain of $\arcsX$ is in fact a fountain by Theorem~\ref{precoveringfountain}.
\end{proof}

\section{Examples}  \label{examplesofarcs}
In this section we are going to describe two special types of torsion theories in $\cal D$, those of t-structures and co-t-structures.

Let us first describe t-structures.

\begin{thm}\label{torsiontheoryconfiguration2}
Let 
$(\cal U, \cal V)$ be a t-structure in $\cal D$, that is, $(\cal U, \cal V)$ is a torsion theory with $\Sigma \cal U \subseteq \cal U$. Suppose $\cal U$ is neither 
zero nor all of $\cal D$, then there is a half line such that the indecomposable objects of $\cal U$ are precisely the objects on the half line and to the left of it, as shown in the following diagram.



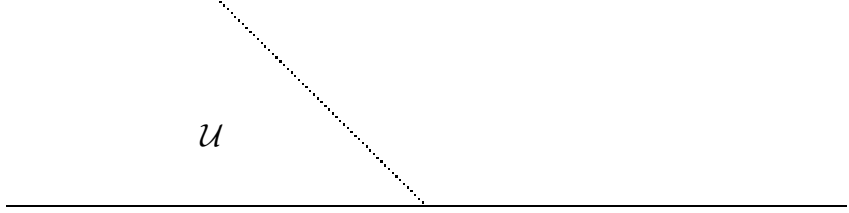
\[
 \xymatrix @-1.0pc @! {
   &  & & & &   & &  & &   & &   & \\
      & &   & & & &   & &   & &   & & \\
    && &   \cal U & &   & &   & &  & &   & & \\
        \ar@{-}[rrrrrrrrrrrr]   & &    & &    &  &  *{} \ar@{.}[uuulll] &  &     & &  &  & \\
              }
\]

\end{thm}



\begin{proof}
Let $\arcsU$ be the corresponding set of arcs for the subcategory $\cal U$. 

(\emph{Step 1})
Consider a horizontal line $y-x=k$ with $k \geq 3$ in the AR quiver of $\cal D$. We claim that, if there are objects from $\cal U$ on this line, then there is a rightmost such object. Namely, suppose not. Then there are objects of $\cal U$ arbitrarily far to the right on $y-x=k$, so all objects on $y-x=k$ are in $\cal U$ because $\Sigma \cal U \subseteq \cal U$. Now, in the following diagram, let $d_1 = (u_1, u_2)$ and $d_2 = (u_1 + 1, u_2 - 1)$ be objects on the lines $y-x = k+1$ and $y-x = k-1$ respectively. Then we have the AR triangle $d_0 \rightarrow  d_1 \oplus d_2 \rightarrow d_0' \rightarrow$, where $d_0 = (u_1, u_2 -1)$ and $d_0' = (u_1+1, u_2)$. Since $d_0$ and $d_0'$ both lie on the line $y-x = k$ which is in $\cal U$, it follows that $d_1$ and $d_2$ are in $\cal U$, since $\cal U$ is closed under extensions and direct summands. Therefore the two neighbouring lines $y-x = k+1$ and $y-x = k-1$ are in $\cal U$. Repeating the argument for other (horizontal) lines, $\cal U$ has to contain all the indecomposable objects of $\cal D$, i.e. $\cal U$ has to be all of $\cal D$. The case where $k = 2$ is similar.

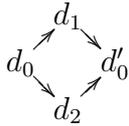
\[
 \xymatrix @-1.5pc @! {
    &  & & &  *{d_1}  \ar[dr] &   & &  & &   & &  &  & \\
     & &   & *{d_0} \ar[ur] \ar[dr] & & *{d_0'}    &   & &   & &   &  &  \\
    && &    &  *{d_2} \ar[ur] &   & &   & &  & &   &   & \\
        \ar@{-}[rrrrrrrrrrrrr]   & &    & &    &  &   &  &     & &  &  & & & & \\
              }
\]

(\emph{Step 2})
Now pick an object $d=(m,n)$ in $\cal U$ and assume that it is the rightmost object of $\cal U$ on the horizontal line $y-x=n-m$. Here we will show that the region $L = \{ (x,y) {\mid} y \leq n, y-x \geq 2\}$ is in $\cal U$.

\begin{itemize}
\item[(i)]Suppose $n-m=2$. In the following diagram, let $u_1 = (m-1, n)$. Then we have the AR triangle $\Sigma d \rightarrow u_1 \rightarrow d \rightarrow$. Since $\Sigma \cal U \subseteq \cal U$, therefore $\Sigma d$ is in $\cal U$. Hence $u_1$ is in $\cal U$, since $\cal U$ is closed under extensions. By applying a (similar) argument on $u_1$ and so on, the half line $y=n$ is in $\cal U$, and so are all the half lines $y=n'$ with $n' \leq n$, since $\Sigma \cal U \subseteq \cal U$. Therefore the region $L$ is in $\cal U$.

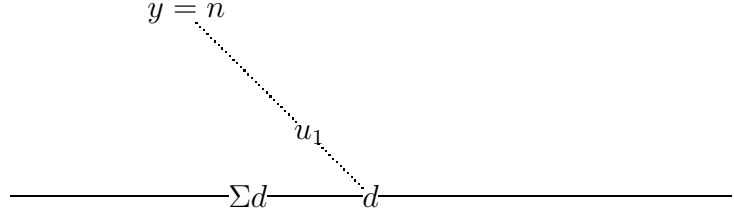
\[
 \xymatrix @-2.5pc @! {
   &  &  &  *{y=n} & &   & &  & &   & &   & \\
      & &   & & & &   & &   & &   & & \\
    && &    & &  *{u_1} \ar@{.}[uull] & &   & &  & &   & & \\
        \ar@{-}[rrrr]   & &    & & *{\Sigma d} \ar@{-}[rr]   &  &  *{d} \ar@{.}[ul] \ar@{-}[rrrrrr] &  &     & &  &  & \\
              }
\]

\item[(ii)]Suppose $n-m > 2$. 
Since $\Sigma \cal U \subseteq \cal U$, therefore $\Sigma d = (m-1,n-1)$ is in $\cal U$. In the following diagram, let $u_1 = (m-1,n)$ and $u_2 = (m, n-1)$. Then we have the AR triangle $ \Sigma d \rightarrow u_1 \oplus u_2 \rightarrow d \rightarrow$. Therefore $u_1$ and $u_2$ are both in $\cal U$, since $\cal U$ is closed under extensions.  By applying a (similar) argument on $u_1$ and $u_2$ and so on (if possible), we can see that the two half lines $t_1$ and $t_2$ are in $\cal U$. Eventually the region labelled $L_0$ (including the two half lines $t_1$ and $t_2$) is in $\cal U$, since $\Sigma \cal U \subseteq \cal U$.

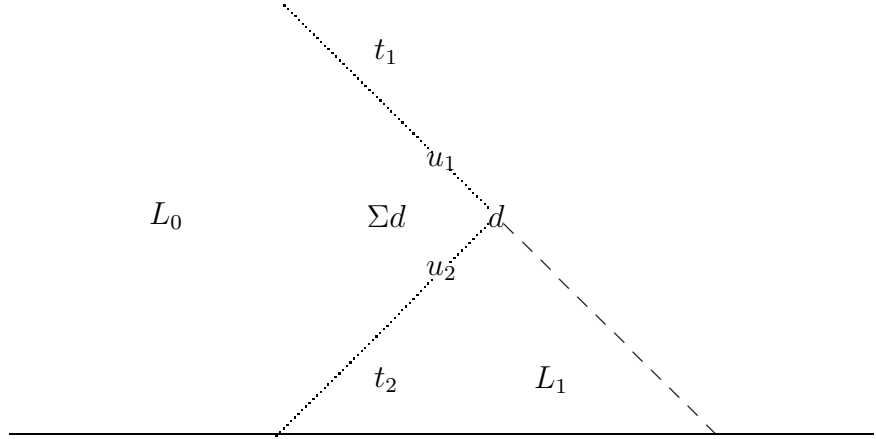
 \[
 \xymatrix @-1.5pc @! {
 & &  &   & &  &  &   & &  & &   & &  & \\
 & &  & &  & &   & *{t_1} &    & &   & &  & &  \\
 & & &   & &   & &  & &    & &   & &   & \\
 &  & & &  & &  &  &  *{u_1} \ar@{.}[uuulll]& &     & &  & &\\
 &  &&  *{L_0} &  &  & &  *{\Sigma d}  &   & *{d}  \ar@{.}[ul]  \ar@{--}[ddddrrrr]  & &   & &   & \\
 &    & & &   & &  & &  *{u_2}   \ar@{.}[ur] & &   & &   & & \\
 &   &&& &   & &   & &   & &  & &   & & \\
 &    && &   & &   & *{t_2} &    & & *{L_1}  &   &   & & \\
        \ar@{-}[rrrrrrrrrrrrrrrr]   & &    & &  & *{} \ar@{.}[uuurrr] &   & &   &&&& & *{} &  &  & \\
              }
\]

Now it remains to show that the little triangular region, $L_1 = \{(x,y) \mid m+1 \leq x \leq n-2, m+3 \leq y \leq n$ and $y-x \geq 2 \}$, is also in $\cal U$.

With the help of the following diagram, let $r_0 = (m-1,m+1)$. Since the arcs $r_0$ and $d=(m,n)$ cross, it follows that $q_0 = (m+1,n)$ is in $\cal U$ since $\arcsU$ satisfies condition (i) by Theorem~\ref{torsiontheoryconfiguration1}. 
Similarly, using $r_1 = (m-1, m+2)$, it follows that $q_1 = (m+2,n)$ is in $\cal U$, and so on (if possible), until all the objects $(m',n)$, with $m+1 \leq m' \leq n-2$, are in $\cal U$. We can repeat this argument, starting with $u_2 =(m, n-1)$, $u_3 =(m, n-2)$ and so on instead (if possible), until the region $L_1$ is in $\cal U$. Therefore $L = L_0 \cup L_1$ is in $\cal U$, and we are done. 

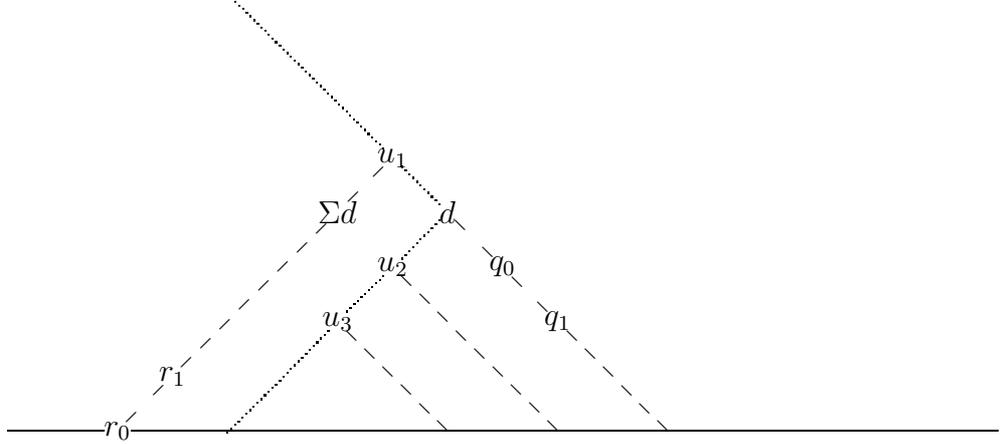
\[
 \xymatrix @-1.5pc @! {
   & &   & &  &  &   & &  & &   & &  & \\
  &  & &  & &   & &    & &   & &  & &  \\
  & &   & &   & &  & &    & &   & &   & \\
    && &  & &  &  &  *{u_1} \ar@{--}[dl] \ar@{.}[uuulll]   & &     & &  & &\\
 &  &   &  & & &  *{\Sigma d}  \ar@{--}[dddlll] &   & *{d} \ar@{.}[ul] \ar@{--}[dr] & &   & &   & \\
     & & &   & &  & &  *{u_2}  \ar@{--}[dddrrr] \ar@{.}[ur] &  & *{q_0}  \ar@{--}[dr] & &   & & \\
    &&& &    & &  *{u_3} \ar@{--}[ddrr] \ar@{.}[ur] & &   & &  *{q_1} \ar@{--}[ddrr] & &   & & \\
    && &  *{r_1} \ar@{--}[dl]  & &   & &   & &  & &   & & \\
        *{}\ar@{-}[rr]   & &   *{r_0} \ar@{-}[rrrrrrrrrrrrrrrr]& & *{} \ar@{.}[uurr]  & *{} & *{} & *{}& *{}&*{}&*{}&*{} &*{}  &*{} &*{}&*{} &  *{}  & *{}& *{} & *{} & *{}\\
              }
\]
\end{itemize}

(\emph{Step 3})
Now suppose there is an indecomposable object $d' = (u,v)$ of $\cal U$ which lies outside the region $L$, i.e. $n < v$. Similar to Step 2, all the half lines $y=v'$ with $v' \leq v$ will be in $\cal U$. Therefore the indecomposable object $ \Sigma^{-1}d = (m+1, n+1)$ will also be in $\cal U$, contradicting that $d$ is the indecomposable object of $\cal U$ which is rightmost on the line  $y-x = n-m$. Therefore the indecomposable objects of $\cal U$ are precisely the objects on the half line $y=n$ and to the left of it, i.e. the region $L$.
\end{proof}

By contrast, there are no non-trivial co-t-structures.

\begin{thm}\label{torsiontheoryconfiguration3}
Let $(\cal U, \cal V)$ be a co-t-structure in $\cal D$, that is, $(\cal U, \cal V)$ is a torsion theory with $\Sigma^{-1} \cal U \subseteq \cal U$. If $\cal U$ is 
non-zero, then $\cal U$ has to be all of $\cal D$. 
\end{thm}

\begin{proof}




Let $\arcsU$ be the corresponding set of arcs for the subcategory $\cal U$. If $\cal U$ is non-zero, let $d = (m,n)$ be an indecomposable object of $\cal U$. 

(\emph{Step 1}) Here we will show that the region $R_0 = \{ (x,y) {\mid} m \leq x, n \leq y$ and $y-x \geq 2\}$ is in $\cal U$.

\begin{itemize}
\item[(i)] Suppose $n-m=2$. 
In the following diagram, let $u_1 = (m, n+1)$. Then we have the AR triangle $d \rightarrow u_1 \rightarrow \Sigma^{-1} d \rightarrow$
. Since $\Sigma^{-1} \cal U \subseteq \cal U$, therefore $\Sigma^{-1} d$ is in $\cal U$. Hence $u_1$ is in $\cal U$, since $\cal U$ is closed under extensions. By applying a similar argument on $u_1$ and so on, the half line $x=m$ is in $\cal U$, and so are all the half lines $x=m'$ with $m \leq m'$, since $\Sigma^{-1} \cal U \subseteq \cal U$. Therefore the region $R_0$ is in $\cal U$.

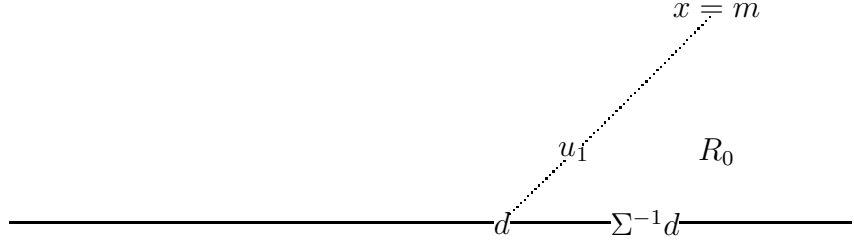
\[
 \xymatrix @-2.5pc @! {
   &  & & & &   & &  & &   &  *{x=m} & & \\
      & &   &  & & &   &  &   & &   & & \\
    && &    & &   &  &   & *{u_1} \ar@{.}[uurr]&   & *{R_0} &   & & \\
      \ar@{-}[rrrrrrr]   & &    & &    &  &    &   *{d} \ar@{.}[ur] \ar@{-}[rr]  &     &  *{\Sigma^{-1} d}\ar@{-}[rrr] &  &  & \\
              }
\]

\item[(ii)] Suppose $n-m>2$. Since $\Sigma^{-1} \cal U \subseteq \cal U$, therefore $\Sigma^{-1} d = (m+1,n+1)$ is in $\cal U$. In the following diagram, let $u_1 = (m,n+1)$ and $u_2 = (m+1, n)$. Then we have the AR triangle $d \rightarrow u_1 \oplus u_2 \rightarrow \Sigma^{-1} d \rightarrow$. Therefore both $u_1$ and $u_2$ are in $\cal U$, since $\cal U$ is closed under extensions. By applying a (similar) argument on $u_1$ and $u_2$ and so on (if possible), we can see that the two half lines $t_1$ and $t_2$ are in $\cal U$. Eventually the region $R_0$ is in $\cal U$, since $\Sigma^{-1} \cal U \subseteq \cal U$. 

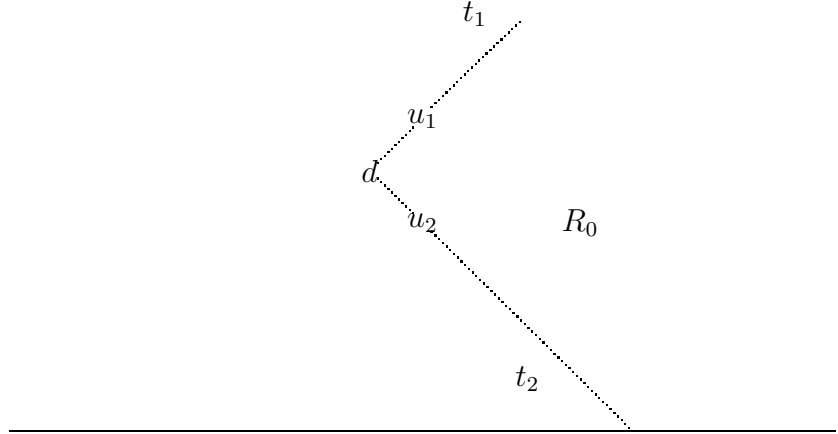
\[
 \xymatrix @-1.5pc @! {
  &  &   & &  &  &   & &  & *{t_1} &   & &  & \\
  &  & &  & &   & &    & &   & &  & &  \\
  & &   & &   & &  & &  *{u_1} \ar@{.}[uurr]  & &   & &   & \\
   & & &  & &  &  &  *{d} \ar@{.}[ur] \ar@{.}[dr] & &     & &  & &\\
   &&   &  &  & &   &   & *{u_2}  \ar@{.}[ddddrrrr] & &   & *{R_0} &   & \\
     & & &   & &  & &   & &   & &   & & \\
    &&& &   & &   & &   & &  & &   & & \\
     && &   & &   & &   & &  & *{t_2} &   & & \\
        \ar@{-}[rrrrrrrrrrrrrrrr]   & &    & &   & &   & &   & && & *{}  & &  &  & \\
              } 
\]
\end{itemize}

(\emph{Step 2}) It follows from Step 1 that $m$ is a right fountain of $\arcsU$. By Theorem~\ref{torsiontheoryconfiguration1}, it is a fountain. 

In particular, choose an indecomposable object $d_0 = (p,m)$ of $\cal U$ with $m-p>2$, as shown in the following diagram. Similar to (ii) in Step 1, the region $R_1 = \{ (x,y) {\mid} p \leq x, m \leq y$ and $y-x \geq 2\} \supset R_0$ is in $\cal U$.

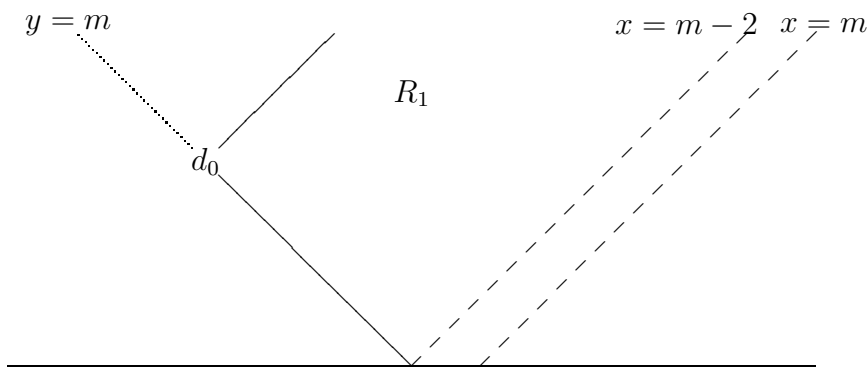
\[
 \xymatrix @-4.3pc @! { 
   &  *{y=m} &  &  & &   & &  & &    &  *{x=m-2} &    & *{x=m} \\
   &   & *{}   &   & &   &  *{R_1} &  & &   &   &   & \\
   &   & &  *{d_0} \ar@{.}[uull]  \ar@{-}[uurr]& &  &   &  & &   &   &   & \\
      & &   &  & *{} & &   &  &   &   &  & & \\
    && &    & & *{}   & &   & &   &  &   & & \\
      \ar@{-}[rrrrrrr]   & &    & &    &  &  *{} \ar@{-}[uuulll] \ar@{--}[uuuuurrrrr] &   *{} \ar@{-}[rrrrr] \ar@{--}[uuuuurrrrr] &     & &  &  & \\
              }
\]

(\emph{Step 3}) Similarly, $p$ is a right fountain of $\arcsU$, and so repeating the argument of Step 2, we see that there is a region $R_2 \supset R_1$ in $\cal U$. Since we can continue in this way indefinitely, $\cal U$ has to contain all the indecomposable objects of $\cal D$, i.e. $\cal U$ has to be all of $\cal D$.

\end{proof}

\medskip

\begin{acknowledgement}
 The author would like to express her thankfulness for the financial assitance from the School of Mathematics and Statistics, Newcastle University, for the Newcastle University International Postgraduate Scholarship (NUIPS), for the Overseas Research Students Awards Scheme (ORSAS) Award, and for the Croucher Foundation Scholarship, which have enabled her study in Newcastle University. The author would also like to express her special gratitude to her supervisor, Peter J\o rgensen, for his generous deliverance of ideas and for the different interpretations he gladly endows on mathematical ideas.
\end{acknowledgement}


\end{document}